\documentclass[11pt]{article}
\usepackage{graphicx}
\usepackage{amsfonts, amsmath, amsthm, amssymb, latexsym, amscd}
\usepackage{pstricks}

\newtheorem{thm}{Theorem}[section]
\newtheorem{exam}[thm]{Example}

\newtheorem{corol}[thm]{Corollary}

\newenvironment{prf}{{\bf \noindent Proof } \rm }{\hfill{$\Box$}}

\title{Inverse single facility location problem in the plane with variable coordinates}

\author{ 
Nazanin Tour-Savadkoohi
	\thanks{Faculty of Mathematical Sciences, Shahrood University of Technology, University Blvd., Shahrood, Iran, e-mail: n.t.savadkoohi@gmail.com.}
, Jafar Fathali
       \thanks{Faculty of Mathematical Sciences, Shahrood University of Technology, University Blvd., Shahrood, Iran, E-mail: fathali@shahroodut.ac.ir.}
}

\begin{document}
\date{}
\maketitle

\begin{abstract}
In traditional facility location problems, a set of points is provided, and the objective is to determine the best location for a new facility based on criteria such as minimizing cost, time, and distances between clients and facilities. Conversely, inverse single facility location problems focus on adjusting the problem's parameters at minimal cost to make a specific point optimal. In this paper, we present an algorithm for the general case of the inverse single facility location problem with variable coordinates in a two-dimensional space. We outline the optimality conditions of this algorithm. Additionally, we examine the specific case namely the inverse minisum single facility location problem and test the algorithm on various instances. The results demonstrate the algorithm's effectiveness in these scenarios.
\end{abstract}


\section{Introduction}
The facility location problem is a well-established topic in operations research and has significant practical implications. Most location models aim to determine the optimal placement of one or more facilities to serve clients, focusing on criteria such as transportation costs and service times. Typically, these models are categorized into three types: continuous, discrete, and network location models.

In continuous models, the goal is to plan facility locations within a two-dimensional space $R^2$. Discrete location models, on the other hand, restrict facility placement to specific predetermined points. Network location models require facilities to be situated on a defined network. Specifically, the continuous version of the single facility location problem seeks to identify a facility's position in a plane that minimizes client servicing costs.

In classical facility location problems, a set of client locations is provided along with assigned weights for each point. A primary objective function in these problems is minimizing the sum of weighted distances between clients and the facility. it is referred to as minisum single facility location problem; it is also known as median problem in location theory literature.

In graph-based versions of facility location models, clients lie on the vertices of a graph, with the objective being to find an optimal facility location within the given network. For further insights into facility location models, readers can refer to \cite{LMV88}.

In certain real-world scenarios, existing facilities may need adjustments rather than relocation. The aim here is to modify problem parameters at minimal cost so that a predetermined point becomes optimal. These scenarios fall under inverse facility location problems. As noted by Heuberger in \cite{H04}, inverse location models become relevant when facilities cannot be relocated without incurring significant costs. In such cases, adjustments can be made by changing the vertices' weight or coordinates with the minimum cost to optimize existing facilities. Inverse location models have applications in banking (see \cite{F22}), hub cargo airport locations (see \cite{AB17}) and train stations (see \cite{AB11}). 


There has been considerable research focused on graph-based inverse facility location problems. For example, Burkard et al. \cite{BPZ04} introduced an $O(n \log n)$ algorithm for addressing the inverse 1-median problem on trees; this was later improved by Galavii \cite{G10}, who achieved linear time complexity for this issue. 
Furthermore, Gassner \cite{Gas10} extends these findings by exploring convex ordered median problems in trees, providing insights into the algorithmic approaches that can be utilized. Other researchers have explored variations of these problems; for instance, Guan and Zhang \cite{GZ12} examined inverse 1-median issue on trees using Hamming norm and Sepasian and Rahbarnia \cite{SR15} tackled the inverse 1-median problem considering varying vertex weights alongside edge lengths on trees. When the underlying network is a block graph, Nguyan \cite{N16-2} proposed a solution algorithm for the inverse 1-median problem with variable vertex weight. The inverse 1-median problem has also been explored in the context of variable edge lengths and weights, as demonstrated by Baroughi et al. \cite{BBG11}. Their findings suggest that the complexity of the problem increases with the introduction of variable parameters, necessitating more sophisticated algorithms to achieve optimal solutions. Similarly, the work of Sepasian and Rahbarnia \cite{SR15} on variable vertex weights and edge reductions further illustrates the adaptability of the inverse 1-median problem to different network configurations. 

The continuous version of the inverse facility location problem has received less attention from researchers than the graph version. Some authors have explored the inverse minisum single facility location problem with variable weights, as seen in papers \cite{BPZ04,BGG10}. Burkard et al. \cite{BPZ04} introduced an algorithm that operates in $O(n\log n)$ time for this problem when the rectilinear norm is used to measure the distances. In contrast, when using the Euclidean norm, Burkard et al. \cite{BGG10} proposed a linear programming model and demonstrated that the unit-cost model can also be solved in $O(n\log n)$ time. Additionally, Baroughi-Bonab et al. \cite{BBA10} examined the inverse minisum single facility location problem with variable coordinates. They established that problems involving rectilinear and Chebyshev norms are NP-hard.

Another continuous inverse location model that has been recently studied is the inverse minisum circle location problem, investigated by Gholami and Fathali \cite{GF22}. In this scenario, the objective is to adjust point weights at minimal cost so that the sum of weighted distances from a given circle's circumference is minimized. The fundamental models for this issue involve infinite constraints. Gholami and Fathali \cite{GF22} proposed a mathematical model with finite constraints; nevertheless, their model is nonlinear and presents challenges for resolution.

Fathali \cite{F22} recently proposed a general row generation algorithm for the continuous inverse location problems with variable weights. The method is starting with solving a linear programming with single constraint and then one constraint is added to the model in each iteration. The added constraint try to change the weights of points such that the given predetermined point become better than the median point respect to the current weights of points. 

In this paper, we develop and modify the method of Fathali \cite{F22} for the variable coordinate case. Since the constrains in the variable case are nonlinear, this case is harder than variable weights version.        


In the subsequent sections, the general formulation of the continuous single facility location problem and its inverse with variable coordinates are discussed in Section \ref{sec2}. Section \ref{sec3} introduces a general algorithm for solving such problems. Following this, Section \ref{secminisum} demonstrates the application of the proposed algorithm to a specific case of continuous inverse location models, namely the inverse minisum single facility location problem. Computational results are presented and analyzed in Section \ref{sec4}. Finally, Section \ref{sec5} provides a summary and concludes the paper.

\section{The facility location problem and its inverse with variable coordinates}\label{sec2}
Let $n$ points $\mathbf{P}_1, \mathbf{P}_2, ..., \mathbf{P}_n$ be given in the plane, where $\mathbf{P}_i=(a_i,b_i)\in \mathbb{R}^2$ for $i=1,...,n$ are the locations of clients and $w_i$ be the weight of point $\mathbf{P}_i$. Set $\mathbf{P}=\{\mathbf{P}_1, \mathbf{P}_2, ..., \mathbf{P}_n\}$ and $\mathbf{w}=(w_1,...,w_n)\geq \mathbf{0}. $
In a single facility location problem, the aim is finding a point $\mathbf{x}=(x,y)$ in the plane such that an objective function $f(\mathbf{x},\mathbf{w},\mathbf{P})$ is minimized, i.e. the following problem should be solved:

\begin{equation*}\label{minimax}
(A_1) \ \min_{\mathbf{x}\in \mathbb{R}^2} f(\mathbf{x},\mathbf{w},\mathbf{P}).
\end{equation*}
The goal of an inverse model is to change some parameters of the problem with minimum cost so that a predetermined point $\bar{\mathbf{x}}$ becomes optimal. Fathali \cite{F22} investigated the case of modifying the weights of points and proposed an efficient algorithm to solve the problem. 

This paper deals with the case of modifying the coordinates of points. Therefore, the problem is to modify the coordinates of existing points $\mathbf{P_i}$, for $i=1,...,n$, to ${\mathbf{\hat{P}_i}}=\mathbf{P_i}+\mathbf{r_i}-\mathbf{s_i}$ with minimum cost, where $\mathbf{r_i}=(r_{i1},r_{i2})$ and $\mathbf{s_i}=(s_{i1},s_{i2})$ are the values of increasing and reducing the coordinate of point $\mathbf{P}_i$.
Let $\mathbf{c^+_i}=(c^+_{i1},c^+_{i2})$ and $\mathbf{c^-_i}=(c^-_{i1},c^-_{i2})$ be the cost vectors of increasing and decreasing per unit of $\mathbf{P_i}$, respectively. Then model of the inverse problem with variable coordinates can be written as follows:
\begin{align}
(A_{2})\ \min\ \ g(\mathbf{r},\mathbf{s})=& \sum_{i=1}^n(c_{i1}^+r_{i1}+c_{i2}^+r_{i2}+c_{i1}^-s_{i1}+c_{i2}^-s_{i2})\label{obj1}\\
&s.t.\notag \\
&f(\bar{\mathbf{x}},\mathbf{w},\mathbf{\hat{P}})\leq f(\mathbf{x},\mathbf{w},\mathbf{\hat{P}}) \ \ \ \ \ \forall\mathbf{x}\in\mathbb{R}^2\label{inf1}\\
&\hat{a_i}=a_i+r_{i1}-s_{i1} \ \ \ \ i=1,\cdots, n,\label{ahat}\\
&\hat{b_i}=b_i+r_{i2}-s_{i2} \ \ \ \ i=1,\cdots, n.\label{bhat}\\
& 0\leq r_{ij}\ \ \ \ ,j=1,2,  \ \ \ \ i=1,\cdots, n,\ \ \ \ \\
& 0\leq s_{ij}\ \ \ \ ,j=1,2,  \ \ \ \ i=1,\cdots, n.\ \ \ \ 
\end{align}
In this model, objective function \eqref{obj1} indicates the minimizing of the cost function. Constraints \eqref{inf1} guarantee that the point $\bar{\mathbf{x}}$ becomes optimal with respect to the coordinates $\mathbf{\hat{P}}$.
Constraints \eqref{ahat} and \eqref{bhat} show the new coordinates of points. 

The model $(A_2)$ is a semi-infinite model and difficult to solve.
In the next section, two efficient algorithms for solving the inverse single facility location problem with variable coordinates are proposed.
 
\section{The proposed algorithms}\label{sec3}
In this section, we propose two algorithms for solving inverse problem with variable coordinates. The basic ideas of our first algorithm is the same as algorithm in \cite{F22}. Where, algorithm starts with finding $\mathbf{x^{(0)}}$, which is the optimal solution of problem $(A_1)$ respect to the current coordinates. Then, in each iteration the algorithm tries to modify the coordinates of points with minimum cost such that the given point $\bar{\mathbf{x}}$ becomes better than all optimal points of problem $(A_1)$ respect to the previous iterations. To achieve this aim, a sub-problem with $t$ constraints should be solved in iteration $t$, i.e. the following problem should be solved, in iteration $t$:   
\begin{align}
	(A_{3}^{(t)})\ \min\ \ g^{(t)}(\mathbf{r},\mathbf{s})= & \sum_{i=1}^n(c_{i1}^+r_{i1}+c_{i2}^+r_{i2}+c_{i1}^-s_{i1}+c_{i2}^-s_{i2})\label{obj}\\
	&s.t.\notag \\
	&f(\bar{\mathbf{x}},\mathbf{w},\mathbf{\hat{P}})\leq f(\mathbf{x^{(k)}},\mathbf{w},\mathbf{\hat{P}}),
	 \ \ \ \ \ for\ k=0,...,t\label{titer}\\
	&\hat{a_i}={a_i}^{(t)}+r_{i1}-s_{i1} \ \ \ \ i=1,\cdots, n,\label{ahat1}\\
	&\hat{b_i}={b_i}^{(t)}+r_{i2}-s_{i2} \ \ \ \ i=1,\cdots, n,\label{bhat1}\\
	& 0\leq r_{ij}\ \ \ \ ,j=1,2,  \ \ \ \ i=1,\cdots, n,\ \ \ \ \\
	& 0\leq s_{ij}\ \ \ \ ,j=1,2,  \ \ \ \ i=1,\cdots, n.\ \ \ \
\end{align}
where $\mathbf{{x^{(t)}}}$ is the solution of problem $(A_1)$ respect to the obtained coordinates in the previous iteration. 
 
Note that in the proposed algorithm, instead of solving the model $(A_2)$ which contains infinite constraints, some sub-problems with finite constraints are solved. Based on above discussion the following algorithm can be used for solving the inverse single facility location problem with variable coordinates.
 
\bigskip
\noindent {\bf Algorithm ISFLP1}\label{al1}.\\
\begin{enumerate}
\item \textbf{Set} $\mathbf{P}^{(0)}=\mathbf{P}$.
\item \textbf{Find} $\mathbf{x^{(0)}}$, the optimal solution of problem $(A_1)$ and set $f_0=f(\mathbf{x}^{(0)},{\mathbf{w}},\mathbf{P}^{(0)})$.
\item \textbf{If} $f(\mathbf{\bar{x}},{\mathbf{w}},\mathbf{P}^{(0)})=f(\mathbf{x}^{(0)},{\mathbf{w}},\mathbf{P}^{(0)})$ then stop (the current solution is optimal). 
\item \textbf{Find} $\mathbf{r_i}^{(0)}$ and $\mathbf{s_i}^{(0)}$, $i=1,...,n$, the optimal solution of problem $(A^{(0)}_3)$.
\item \textbf{For} $i=1,...,n$, set $\mathbf{P_i^{(1)}}=\mathbf{P_i^{(0)}}+{\mathbf{r_i}^{(0)}}-\mathbf{s_i}^{(0)}$.
\item \textbf{Find} $\mathbf{x^{(1)}}$, the optimal solution of the following sub-problem:
$$\min_{\mathbf{x}\in \mathbb{R}^2} f({\mathbf{x}},{\mathbf{w}},\mathbf{P}^{(1)}),$$
 and set $f_1=f(\mathbf{x}^{(1)},{\mathbf{w}},\mathbf{P}^{(1)})$.
\item \textbf{Set} $k=1$.
\item \textbf{While} $\frac{|f(\mathbf{\bar{x}},\mathbf{w},\mathbf{\hat{P}})-f(\mathbf{x}^{(k)},\mathbf{w},\mathbf{\hat{P}})|}{|f(\mathbf{\bar{x}},\mathbf{w},\mathbf{\hat{P}})|}> \epsilon$ do
\begin{enumerate}
\item \textbf{Find} $\mathbf{r_i}^{(k)}$ and $\mathbf{s_i}^{(k)}$, $i=1,...,n$, the optimal solution of problem $(A^{(k)}_3)$.\\

\item \textbf{For} $i=1,...,n$, set $\mathbf{P_i}^{(k+1)}=\mathbf{P_i}^{(k)}+{\mathbf{r_i}^{(k)}}-\mathbf{s_i}^{(k)}$.
\item \textbf{Set} $k=k+1$.
\item \textbf{Find} $\mathbf{x^{(k)}}$, the optimal solution of the following sub-problem:
$$\min_{\mathbf{x}\in \mathbb{R}^2} f({\mathbf{x}},{\mathbf{w}},\mathbf{P}^{(k)}),$$ 
and $f_k=f(\mathbf{x}^{(k)},{\mathbf{w}},\mathbf{P}^{(k)})$.
\end{enumerate}
\textbf{End while}
\item \textbf{For} $i=1,...,n$ set
$$(r_{i1}^*,r_{i2}^*,s_{i1}^*,s_{i2}^*)=(\sum_{i=1}^{n}{r_{i1}^{(k-1)}},\sum_{i=1}^{n}{r_{i2}^{(k-1)}} ,\sum_{i=1}^{n}{s_{i1}^{(k-1)}},\sum_{i=1}^{n}{s_{i2}^{(k-1)}}), $$
and $C^*=g^{(k-1)}(\mathbf{r^*},\mathbf{s^*})$.
\end{enumerate}
\textbf{End of algorithm}
\bigskip

In the second algorithm, which we call it ISFLP2, to improve the time complexity, we replace the constraints \eqref{titer} by the last single constraint. This algorithm is the same as ISFLP1, just in step 8.a, the following problem should be solved:  
\begin{align}
	(A_{4}^{(t)})\ \min\ \ g^{(t)}(\mathbf{r},\mathbf{s})= & \sum_{i=1}^n(c_{i1}^+r_{i1}+c_{i2}^+r_{i2}+c_{i1}^-s_{i1}+c_{i2}^-s_{i2})\label{obj}\\
	&s.t.\notag \\
	&f(\bar{\mathbf{x}},\mathbf{w},\mathbf{\hat{P}})\leq f(\mathbf{x^{(t)}},\mathbf{w},\mathbf{\hat{P}}),\label{siter}
	\ \ \ \ \ \\
	&\hat{a_i}={a_i}^{(t)}+r_{i1}-s_{i1} \ \ \ \ i=1,\cdots,n,\\
	&\hat{b_i}={b_i}^{(t)}+r_{i2}-s_{i2} \ \ \ \ i=1,\cdots, n,\\
	&0\leq {r}_{ij}\ \ \ \ ,j=1,2,\  \ \ \ \ i=1,\cdots, n, \\
	&0\leq {s}_{ij}\ \ \ \ ,j=1,2,\ \ \ \ \ i=1,\cdots, n.
\end{align}

The following theorem shows that if in an iteration the value of objective function in sub-problem is equal to $f(\mathbf{\bar{x}},{\mathbf{w}},\mathbf{P}^{(t)})$ then an optimal solution is found. 

\begin{thm}\label{stopfx}
	If in iteration $t$ of Algorithms ISFLP1 and ISFLP2, $f(\mathbf{\bar{x}},{\mathbf{w}},\mathbf{P}^{(t)})=f(\mathbf{x}^{(t)},{\mathbf{w}},\mathbf{P}^{(t)})$, then $\mathbf{P}^{(t)}$ is an optimal solution of $(A_2)$.
\end{thm}
\begin{prf}
	Note that $\mathbf{x}^{(t)}$ is the optimal solution of the sub-problem, thus
	$$f(\mathbf{x}^{(t)},{\mathbf{w}},\mathbf{P}^{(t)})=\min_{\mathbf{x}\in \mathbb{R}^2} f(\mathbf{x},{\mathbf{w}},\mathbf{P}^{(t)}).$$ 
	Since  $f(\mathbf{\bar{x}},{\mathbf{w}},\mathbf{P}^{(t)})=f(\mathbf{x}^{(t)},{\mathbf{w}},\mathbf{P}^{(t)})$, therefore 
	$$f(\bar{\mathbf{x}},{\mathbf{w}},\mathbf{P}^{(t)})\leq f(\mathbf{x},\mathbf{w},\mathbf{P}^{(t)})\ \ \forall\mathbf{x}\in\mathbb{R}^2,$$ 
	so $\mathbf{P}^{(t)}$ is a feasible solution for problem $(A_2)$.\\
  Now let $(\mathbf{P'},\mathbf{r'},\mathbf{s'})$ be a feasible solution of $(A_2)$. Then 
 $$f(\bar{\mathbf{x}},{\mathbf{w}},\mathbf{P'})\leq f(\mathbf{x},{\mathbf{w}},\mathbf{P'}) \ \ \forall\mathbf{x}\in\mathbb{R}^2.$$
 This implies that $(\mathbf{P'},\mathbf{r'},\mathbf{s'})$ is feasible for $(A_3^{(t)})$, too. Since  $(\mathbf{P}^{(t)},\mathbf{r^{(t)}},\mathbf{s^{(t)}})$ is the optimal solution of $(A^{(t)}_3)$ then 
\[ \sum_{i=1}^n(c_{i1}^+r_{i1}+c_{i2}^+r_{i2}+c_{i1}^-s_{i1}+c_{i2}^-s_{i2})
\leq \sum_{i=1}^n(c_{i1}^+r'_{i1}+c_{i2}^+r'_{i2}+c_{i1}^-s'_{i1}+c_{i2}^-s'_{i2}).\]
 It means $(\mathbf{P}^{(t)},\mathbf{r^{(t)}},\mathbf{s^{(t)}})$ is the optimal solution of $(A_2)$.\\
 The proof for Algorithm ISFLP2 is the same.
\end{prf}   

Note that if $\mathbf{\bar{x}}= \mathbf{x}^{(t)}$ then $f(\mathbf{\bar{x}},{\mathbf{w}},\mathbf{P}^{(t)})=f(\mathbf{x}^{(t)},{\mathbf{w}},\mathbf{P}^{(t)})$. Thus the following property can be concluded.  

\begin{corol}\label{stopx}
	If in iteration $t$ of Algorithms ISFLP1 and ISFLP2, $\mathbf{\bar{x}}= \mathbf{x}^{(t)}$, then $\mathbf{P}^{(t)}$ is optimal for $(A_2)$ .
\end{corol}

The stopping condition of the algorithms are based on Theorem \ref{stopfx}. This stopping condition can be replaced by another one which is the case that the coordinates of points do not change.

\begin{thm}\label{stopw}
Let $t$ be an iteration of Algorithms ISFLP1 and ISFLP2 such that $\mathbf{P}^{(t+1)}=\mathbf{P}^{(t)}$. Then $\mathbf{{P}}^{(t)}$ is the optimal solution of $(A_2)$.
\end{thm}
\begin{prf}
In iteration $t$ of Algorithm ISFLP1, the problem $(A_3^{(t)})$ is solved,
where $\mathbf{x^{(t)}}$ is the solution of $t$-th sub-problem. Thus
$$f({\mathbf{x}^{(t)}},{\mathbf{w}},\mathbf{P}^{(t)})=\min_{\mathbf{x}\in \mathbb{R}^2}  f(\mathbf{x},{\mathbf{w}},\mathbf{P}^{(t)}).$$ 
By the algorithm $(\mathbf{P}^{(t+1)},\mathbf{r}^{(t)},\mathbf{s}^{(t)})$ is the optimal solution of $(A_3^{(t)})$. Then 
$$f(\bar{\mathbf{x}},{\mathbf{w}},\mathbf{P}^{(t+1)})\leq f(\mathbf{x}^{(t)},{\mathbf{w}},\mathbf{P}^{(t+1)}),$$
$$\mathbf{P}^{(t+1)}=\mathbf{P}^{(t)}+\mathbf{r}^{(t)}-\mathbf{s}^{(t)}. $$
Since $\mathbf{P}^{(t)}=\mathbf{P}^{(t+1)}$ then  
$$f(\mathbf{x}^{(t)},{\mathbf{w}},\mathbf{P}^{(t+1)})=\min_{\mathbf{x}\in \mathbb{R}^2}  f(\mathbf{x},{\mathbf{w}},\mathbf{P}^{(t+1)}).$$ 
Therefore,
$$f(\bar{\mathbf{x}},{\mathbf{w}},\mathbf{P}^{(t+1)})\leq f(\mathbf{x},{\mathbf{w}},\mathbf{P}^{(t+1)})\ \ \forall\mathbf{x}\in\mathbb{R}^2,$$ 
It means $(\mathbf{P}^{(t+1)},\mathbf{r}^{(t)},\mathbf{s}^{(t)})$ is a feasible solution of $(A_2)$.
The proof of optimality is the same as proof of Theorem \ref{stopfx}. With the same manner the theorem can be solved for Algorithm ISFLP2.
\end{prf}

Note that, the feasibility is improved in each iteration of algorithm. Therefore,  the algorithm terminated in a near feasible solution with the best objective function.  

The time complexity of the algorithms depend on complexity of sub-problem $(A_1)$ and problems $(A_3^{(t)})$ and $(A_4^{(t)})$, in iteration $t$.   In the next sections we explain the presented algorithms for solving the inverse minisum single facility location problem.

\section{Inverse minisum single facility location problem}\label{secminisum}
In this section, we consider the inverse Minisum Single Facility Location Problem (MSFLP). The objective function of MSFLP, which is also called the Fermat-Weber problem, is minimizing the weighted sum of distances between existing points and the new facility, i.e. 
\begin{equation*} 
(A_{ms}) \ \min_{\mathbf{x}\in \mathbb{R}^2} F(\mathbf{x},\mathbf{w},\mathbf{P})=\sum_{i=1}^n w_id(\mathbf{x},\mathbf{P}_{i}),
\end{equation*}
where $d(\mathbf{x},\mathbf{P}_i)$ is the distance between the points $\mathbf{x}$ and $\mathbf{P}_i$. 

Thus, the inverse of MSFLP with variable coordinates can be written as follows:
\begin{align}
(A_{5})\ \min\ \ g(\mathbf{r},\mathbf{s})= & \sum_{i=1}^n(c_{i1}^+r_{i1}+c_{i2}^+r_{i2}+c_{i1}^-s_{i1}+c_{i2}^-s_{i2})\label{obj}\\	
&s.t.\notag \\
&\sum_{i=1}^{n}{w_i} d(\bar{\mathbf{{x}}},\mathbf{\hat{P}}_i)\leq \sum_{i=1}^{n}{w_i} d(\mathbf{{x}},\mathbf{\hat{P}}_i) \ \ \ \ \ \forall\mathbf{x}\in\mathbb{R}^2\label{infms},\\
	&\hat{a_i}=a_i+r_{i1}-s_{i1} \ \ \ \ i=1,\cdots,n,\\
&\hat{b_i}=b_i+r_{i2}-s_{i2} \ \ \ \ i=1,\cdots, n, \\
&0\leq r_{ij}\ \ \ \ ,j=1,2,\ \ \ \ i=1,\cdots, n,\\
&0\leq s_{ij}\ \ \ \ ,j=1,2,\ \ \ \ i=1,\cdots, n.
\end{align} 
The paper of Baroughi-Bonab et al. \cite{BBA10} is the only paper devoted on inverse of MSFLP with variable coordinates. They showed the problem under rectilinear and Chebyshev norms is NP-hard. In case that the squared Euclidean norm is used, they proposed a linear time algorithm. 

Using Algorithms ISFLP1 and ISFLP2, in iteration $k$, the following sub-problem is solved:
\begin{equation}\label{minisum-sub}
\min_{\mathbf{x}\in \mathbb{R}^2} \sum_{i=1}^n w_id(\mathbf{x},\mathbf{P}_i^{(k)}).
\end{equation}

This sub-problem can be solved by the methods presented in \cite{LMV88}.
Specially, for the squared Euclidean norm, the optimal solution $x^*,y^*$ is obtained by the following formula, which is called center of gravity method.  

\[x^*=\frac{\sum_{i=1}^{n}w_i a_i}{\sum_{i=1}^{n}w_i}, \ \  \ y^*=\frac{\sum_{i=1}^{n}w_i b_i}{\sum_{i=1}^{n}w_i}. \]

Moreover, in iteration $k=t$, the constraints \eqref{titer} of problem ${A_3^{(t)}}$ and \eqref{siter} of problem ${A_4^{(t)}}$ are as follows, respectively.
\begin{equation}\label{minisum-const1}
	\sum_{i=1}^n {w_{i}}(d(\mathbf{\bar{x}},\mathbf{\hat{P}_i})-d(\mathbf{x}^{(j)},\mathbf{\hat{P}_i}))\leq 0 \ \ \ j=0,...t.
\end{equation} 

\begin{equation}\label{minisum-const2}
\sum_{i=1}^n {w_{i}}(d(\mathbf{\bar{x}},\mathbf{\hat{P}_i})-d(\mathbf{x}^{(t)},\mathbf{\hat{P}_i}))\leq 0.
\end{equation} 
We consider the cases that distances are measured by squared Euclidean, Euclidean and rectilinear norms. In these cases, constraint \eqref{minisum-const2} is converted to the following constraints, respectively:
\begin{equation}\label{minisum-constse}
	\sum_{i=1}^n {w_{i}}((\bar{x}-\hat{a}_i)^2+(\bar{y}-\hat{b}_i)^2-
{({x^{(t)}}-\hat{a}_i)^2-({y^{(t)}}-\hat{b}_i)^2})\leq 0,
\end{equation}  

\begin{equation}\label{minisum-conste}
	\sum_{i=1}^n {w_{i}}(\sqrt{(\bar{x}-\hat{a}_i)^2+(\bar{y}-\hat{b}_i)^2}-
	\sqrt{({x^{(t)}}-\hat{a}_i)^2+({y^{(t)}}-\hat{b}_i)^2}\ )\leq 0,
\end{equation}  

\begin{equation}\label{minisum-constr}
	\sum_{i=1}^n {w_{i}}(|\bar{x}-\hat{a}_i|+|\bar{y}-\hat{b}_i|-
	|{x^{(t)}}-\hat{a}_i|-|y^{(t)}-\hat{b}_i|)\leq 0.
\end{equation}  

Constraint \eqref{minisum-constse} can be simplified to the following linear constraint.
\begin{equation}\label{minisum-constse1}
	\sum_{i=1}^n {w_{i}}(\bar{x}^2-2\bar{x}\hat{a}_i+\bar{y}^2-2\bar{y}\hat{b}_i-
	{(x^{(t)})}^2+2x^{(t)}\hat{a}_i-{(y^{(t)})}^2+2y^{(t)}\hat{b}_i)\leq 0.
\end{equation} 

Thus in the case squared Euclidean, the models ${A_3^{(t)}}$ and ${A_4^{(t)}}$ are linear programming. Therefore, in this case the time complexity of each iteration of both proposed algorithms are polynomial.

\begin{exam}\label{ex1}
Consider Table \ref{inst1} which contains the data of an instance of the inverse minisum location problem. Let the given point be $\bar{\mathbf{x}}=(0,1)$.
We examined our algorithms for finding the solution of this instance for squared Euclidean, Euclidean and rectilinear norms. The results are shown in Tables \ref{result-SE-4points} to \ref{result-L1-4points}. In these tables the column with heading $\mathbf{x}^{(k)}$ indicates the solution of sub-problem (\ref{minisum-sub}) in iteration $k$ with respect to $\mathbf{P}^{(k)}$.  $F_{k}=F({\mathbf{x}}^{(k)},\mathbf{w},{\mathbf{P}}^{(k)})$ and $F_B=F({\mathbf{x}}_{B},\mathbf{w},{\mathbf{P}}_{B})$ indicate the values of objective function of minisum location problem in points $\mathbf{{x}^{(k)}}$ and $\mathbf{{x}_{B}}$, respectively.

\begin{table}\caption{ The parameters of the considered instance in Example \ref{ex1}}\label{inst1}.
\centering
\begin{tabular}{ccccccc}
\hline
i& $A_i=(a_i,b_i)$& $w_i$& $c_{i1}^{+}$&$c_{i2}^{+}$&$c_{i1}^{-}$&$c_{i2}^{-}$\\
\hline
1& $(1,0)$&6&$\sqrt{2}$&5&$1$&$6$\\
2& $(-5,3)$&3&7&3&4&2\\ 
3& $ (7,2)$&1&1&2&4&4\\
4& $(0,-0.5)$&2&2&1&4&1\\
\hline
\end{tabular}
\end{table}

\begin{itemize}
	\item For squared Euclidean norm, the iteration results of Algorithms ISFLP1 and ISFLP2 are the same. They are reported in Table \ref{result-SE-4points}. The results show $\mathbf{x}^{(10)}$ converges to the given point $\mathbf{\bar{x}}$. The algorithms find the optimal solution with tolerance $||\bar{\mathbf{x}}-{\mathbf{x}}^{(k)}||\leq\epsilon$, where $\epsilon=
	0.0001$. The cost of modifying $\mathbf{P}^{(0)}$ to $\mathbf{P}^{(10)}$ is $C^*=1.4709$.  Table \ref{result-B-4points} contains the results of  algorithm of Baroughi-Bonab et al. \cite{BBA10} with $C^*= ?? $. 
	
	\begin{table}\caption{ The iteration results of Algorithms  ISFLP1 and ISFLP2 run on instance of Example \ref{ex1} for squared Euclidean.}\label{result-SE-4points}
		\centering
		\tiny 
		\begin{tabular}{ccccccc}
			\hline
			iteration& $\mathbf{x}^{(k)}=(x^{(k)},y^{(k)})$& $\mathbf{P}^{(k)}_1$& $\mathbf{P}^{(k)}_2$&$\mathbf{P}^{(k)}_3$&$\mathbf{P}^{(k)}_4$&$||\bar{\mathbf{x}}-\mathbf{x}^{(k)}||$\\
			\hline
			
			0&(-0.1667 , 0.8333)
			&(1.0000,0.0000)&(-5.0000,3.0000)&(7.0000,2.0000)
			&(0.0000,-0.5000)
			&0.23570\\
			1& (0.0000 , 0.8333)
			&(1.3333,0.0000)&(-5.0000,3.0000)&(7.0000,2.0000)
			&(0.0000,-0.5000)
			&0.16667\\
			2& (0.0000 , 0.9167)
			&(1.3333,0.0000)&(-5.0000,3.0000)&(7.0000,2.0000)
			&(0.0000,0.0000)
			&0.08333\\ 
			3& (0.0000 , 0.9583)&(1.3333,0.0000)&(-5.0000,3.0000)&(7.0000,2.0000)
			&(0.0000,0.2500)
			&0.04167\\
			4& (0.0000 , 0.9792)&(1.3333,0.0000)&(-5.0000,3.0000)&(7.0000,2.0000)
			&(0.0000,0.3750)
			&0.02083\\
			5& (0.0000 , 0.9896)
			&(1.3333,0.0000)&(-5.0000,3.0000)&(7.0000,2.0000)
			&(0.0000,0.4375)
			&0.01042\\
			6& (0.0000 , 0.9948)
			&(1.3333,0.0000)&(-5.0000,3.0000)&(7.0000,2.0000)
			&(0.0000,0.4687)
			&0.00521\\ 
			7&  (0.0000 , 0.9974)&(1.3333,0.0000)&(-5.0000,3.0000)&(7.0000,2.0000)
			&(0.0000,0.4844)
			&0.00260\\
			8& (0.0000 , 0.9987)
			&(1.3333,0.0000)&(-5.0000,3.0000)&(7.0000,2.0000)
			&(0.0000,0.4922)
			&0.00130\\
			9& (0.0000 , 0.9993)
			&(1.3333,0.0000)&(-5.0000,3.0000)&(7.0000,2.0000)
			&(0.0000,0.4961)
			&0.00130\\
			10&(0.0000 , 0.9997) &(1.3333,0.0000)&(-5.0000,3.0000)&(7.0000,2.0000)
			&(0.0000,0.4980)
			&0.00032\\ 
			\hline
			$F_{10}$&154.1676&&&&&\\
			\hline
			&CPU(in s) of ISFLP1&2.1392&&&&\\
			&CPU(in s) of ISFLP2&1.0366&&&&\\
			\hline
		\end{tabular}
	\end{table}

	\begin{table}\caption{The results of method of Baroughi-Bonab et al. \cite{BBA10} on instance of Example \ref{ex1}.}\label{result-B-4points}
		\centering
		\tiny 
		\begin{tabular}{cccccc}
			\hline
			 $\mathbf{x}_{B}=(x_{B},y_{B})$& $\mathbf{P}_{1B}$& $\mathbf{P}_{2B}$&$\mathbf{P}_{3B}$&$\mathbf{P}_{4B}$&$||\bar{\mathbf{x}}-\mathbf{x}_{B}||$\\
			\hline
			( 0.0000, 1.0000)&
		(1.3333,0.0000)&(-5.0000,3.0000)&(7.0000,2.0000)&
		(0.0000,0.5000)&
		0.00000\\
			\hline
		$F_B$&154.1667&&&&\\
			\hline
		CPU(in s)&0.0426&&&&\\
			
			\hline
		\end{tabular}
	\end{table}

	\item  For Euclidean norm, the same results were obtained from both Algorithms ISFLP1 and ISFLP2 as shown in Table \ref{result-E-4points}. The algorithms find the optimal solution with tolerance $\dfrac{|\bar{F}-F_k|}{|\bar{F}|}\leq\epsilon$, where $\epsilon=0.01$. In this case, the cost of modifying $\mathbf{P}^{(0)}$ to $\mathbf{P}^{(9)}$ is $C^*=5.7684$.

\begin{table}\caption{ The iteration results of Algorithms  ISFLP1 and ISFLP2 run on instance of Example \ref{ex1} for Euclidean norm.}\label{result-E-4points}
	\centering
	\tiny 
	\begin{tabular}{ccccccc}
		\hline
		iteration& $\mathbf{x}^{(k)}=(x^{(k)},y^{(k)})$& $\mathbf{P}^{(k)}_1$& $\mathbf{P}^{(k)}_2$&$\mathbf{P}^{(k)}_3$&$\mathbf{P}^{(k)}_4$&$\dfrac{|\bar{F}-F_k|}{|\bar{F}|}$\\
		\hline
		
		0&(0.9768,0.0049)
		&(1.0000,0.0000)&(-5.0000,3.0000)&(7.0000,2.0000)
		&(0.0000,-0.5000)
		&0.1720\\
		1& (0.3059,0.0003)
		&(0.3219,0.0000)&(-5.0000,3.0000)&(7.0000,2.0000)
		&(0.0000,-0.5000)
		&0.1841\\
		2&(-1.2868,0.0075)
		&(-1.2895,0.0000)&(-5.0000,3.0000)&(7.0000,2.0000)
		&(0.0000,-0.5000)
		&0.2880\\ 
		3&(-0.1728,0.0002)
		&(-0.1678,0.0000)&(-5.0000,3.0000)&(7.0000,2.0000)
		&(0.0000,-0.5000)
		&0.2082\\
		4&(-0.1753,0.5694)
		&(-0.1678,0.5717)&(-5.0000,3.0000)&(7.0000,2.0000)
		&(0.0000,-0.5000)
		&0.1121\\
		5&(0.0688,0.7719)
		&(0.0793,0.7754)&(-5.0000,3.0000)&(7.0000,2.0000)
		&(0.0000,-0.5000)
		&0.0503\\
		6&(-0.3684,0.7729)
		&(-0.3621,0.7754)&(-5.0000,3.0000)&(7.0000,2.0000)
		&(0.0000,-0.5000)
		&0.1114\\ 
		7&(-0.03671,0.7720)
		&(-0.0273,0.7754)&(-5.0000,3.0000)&(7.0000,2.0000)
		&(0.0000,-0.5000)
		&0.0549\\
		8&(-0.0369,0.9022)
		&(-0.0273,0.9059)&(-5.0000,3.0000)&(7.0000,2.0000)
		&(0.0000,-0.5000)
		&0.0246\\
		9&(-0.0040,0.9509)
		&(0.0060,0.9549)&(-5.0000,3.0000)&(7.0000,2.0000)
		&(0.0000,-0.5000)
		&0.0095\\
		\hline
		$F_9$&26.2478&&&&&\\
		\hline
		&CPU(in s) of ISFLP1&24.2648&&&&\\
		&CPU(in s) of ISFLP2&11.2763&&&&\\
		\hline
	\end{tabular}
\end{table}

\item  In rectilinear norm case, the results for Algorithms ISFLP1 and ISFLP2 are reported in Table \ref{result-L1-4points}. The both algorithms find the optimal solution with tolerance $\dfrac{|\bar{F}-F_k|}{|\bar{F}|}\leq\epsilon$, where $\epsilon=0.01$.

\begin{table}\caption{ The iteration results of Algorithms  ISFLP1 and ISFLP2 run on instance of Example \ref{ex1} for rectilinear norm.}\label{result-L1-4points}
	\centering
	\tiny 
	\begin{tabular}{ccccccc}
		\hline
		iteration& $\mathbf{x}^{(k)}=(x^{(k)},y^{(k)})$& $\mathbf{P}^{(k)}_1$& $\mathbf{P}^{(k)}_2$&$\mathbf{P}^{(k)}_3$&$\mathbf{P}^{(k)}_4$&$\dfrac{|\bar{F}-F_k|}{|\bar{F}|}$\\
			\hline
		&&&ISFLP1 &&&\\
		\hline
		
		0&(1.0000,0.0000)
		&(1.0000,0.0000)&(-5.0000,3.0000)&(7.0000,2.0000)
		&(0.0000,-0.5000)
		&0.1364\\
		1&(0.5000,0.0000)
		&(0.5000,0.0000)&(-5.0000,3.0000)&(7.0000,2.0000)
		&(0.0000,-0.5000)
		&0.1220\\
		2&(0.0833,0.0000)
		&(0.0833,0.0000)&(-5.0000,3.0000)&(7.0000,2.0000)
		&(0.0000,-0.5000)
		&0.1082\\ 
		3&(0.0000,0.2640)
		&(0.0000,0.2640)&(-5.0000,3.0000)&(7.0000,2.0000)
		&(0.0000,-0.5000)
		&0.0808\\
		4&(0.0000,0.5093)
		&(0.0000,0.5093)&(-5.0000,3.0000)&(7.0000,2.0001)
		&(0.0000,-0.5000)
		&0.0562\\
		5&(0.0000,0.6729)
		&(0.0000,0.6729)&(-5.0000,3.0000)&(7.0000,2.0001)
		&(0.0000,-0.5000)
		&0.0385\\
		6&(0.0000,0.7819)
		&(0.0000,0.7819)&(-5.0000,3.0000)&(7.0000,2.0001)
		&(0.0000,-0.5000)
		&0.0262\\ 
		7&(0.0000,0.8546)
		&(0.0000,0.8546)&(-5.0000,3.0000)&(7.0000,2.0001)
		&(0.0000,-0.5000)
		&0.0177\\
		8&(0.0000,0.9031)
		&(0.0000,0.9031)&(-5.0000,3.0000)&(7.0000,2.0001)
		&(0.0000,-0.5000)
		&0.0119\\
		9&(0.0000,0.9354)
		&(0.0000,0.9354)&(-5.0000,3.0000)&(7.0000,2.0001)
		&(0.0000,-0.5000)
		&0.0080\\
		\hline
		$F_9$&32.1292&&&&&\\
		\hline
		$C^*$&5.6772&&&&&\\
		\hline
		&CPU(in s) of ISFLP1&11.1709
		&&&&\\		
		\hline
		&&&ISFLP2 &&&\\
		\hline
		0&(1.0000,0.0000)
		&(1.0000,0.0000)&(-5.0000,3.0000)&(7.0000,2.0000)
		&(0.0000,-0.5000)
		&0.1364\\
		1&(0.5000,0.0000)
		&(0.5000,0.0000)&(-5.0000,3.0000)&(7.0000,2.0000)
		&(0.0000,-0.5000)
		&0.1220\\
		2&(0.0833,0.0000)
		&(0.0833,0.0000)&(-5.0000,3.0000)&(7.0000,2.0000)
		&(0.0000,-0.5000)
		&0.1082\\ 
		3&(0.000 0,0.2647)
		&(0.0000,0.2647)&(-5.0004,2.9993)&(7.0000,2.0010)
		&(0.0000,-0.5000)
		&0.0808\\
		4&(0.0000,0.5098)
		&(0.0000,0.5098)&(-5.0004,2.9993)&(7.0000,2.0010)
		&(0.0000,-0.5000)
		&0.0561\\
		5&(0.0000,0.6732)
		&(0.0000,0.6732)&(-5.0004,2.9993)&(7.0000,2.0010)
		&(0.0000,-0.5000)
		&0.0385\\
		6&(0.0000,0.7821)
		&(0.0000,0.7821)&(-5.0004,2.9993)&(7.0000,2.0010)
		&(0.0000,-0.5000)
		&0.0262\\ 
		7&(0.0000,0.8548)
		&(0.0000,0.8548)&(-5.0004,2.9993)&(7.0000,2.0010)
		&(0.0000,-0.5000)
		&0.0117\\
		8&(0.0000,0.9032)
		&(0.0000,0.9032)&(-5.0004,2.9993)&(7.0000,2.0010)
		&(0.0000,-0.5000)
		&0.0119\\
		9&(0.0000,0.9354)
		&(0.0000,0.9354)&(-5.0004,2.9993)&(7.0000,2.0010)
		&(0.0000,-0.5000)
		&0.0080\\
		\hline
		$F_9$&32.1294&&&&&\\
		\hline
	$C^*$&5.6823&&&&&\\
		\hline
		&CPU(in s) of ISFLP2&14.9403&&&&\\
		\hline
		
	\end{tabular}
\end{table}
	\end{itemize}
\end{exam} 

\section{Computational results}\label{sec4}
In this section, we evaluate the performance of our proposed algorithms using several test problems. All experiments were conducted on a PC equipped with an Intel Core i3 processor, 4 GB of RAM, and a CPU clocked at 2.7 GHz. The algorithms designed for solving inverse minisum single facility location problem were tested on three distinct problems, each involving different sets of given points.

For the first test problem, we analyze an instance comprising 18 existing points. The coordinates of these points, along with their associated data, are provided in Table \ref{inst-data-18}. The point coordinates were sourced from \cite{F09}.

\begin{table}[h]\centering
	\small
	\begin{tabular}{cc|cc|c}\hline
		$(a_i,b_i,w_i,c_{i1}^+,c_{i2}^+,c_{i1}^-,c_{i2}^-)$  & &$(a_i,b_i,w_i,c_{i1}^+,c_{i2}^+,c_{i1}^-,c_{i2}^-)$& &$(a_i,b_i,w_i,c_{i1}^+,c_{i2}^+,c_{i1}^-,c_{i2}^-)$\\\hline
		$(1,2,3,1.5,2,1,2)$& &$(4,4,1,2,1.7,4,9)$& &$(7,1,2,4,2,2,4)$\\
		$(1,3,2,3,1,3,4)$& &$(4,9,2,2.5,3,1,1)$& &$(7,2,3,3,1,4,3)$\\
		$(2,5,1,6,4,5,6)$& &$(5,3,2,1,5,2,3)$& &$(8,5,1,2,3,5,5)$\\
		$(3,6,3,3,1,7,2)$& &$(5,5,1,2.5,1,4,2)$& &$(8,8,3,3,4,3,7)$\\
		$(4,8,2,1,2,3,1)$& &$(6,6,3,3,2,5,6)$& &$(9,7,3,1,2.5,8,1)$\\
		$(4,1,3,5,4,5,6)$& &$(6,3,3,4,3,4,2)$& &$(9,6,2,1.5,4,2,3)$\\\hline
	\end{tabular}
	\caption{Data for the problem with 18 given points.}\label{inst-data-18}
\end{table}
We also evaluated the proposed algorithms on two additional instances: Ruspini 75 and Bongartz 287, which are available in Beasley \cite{B90}. The weights and costs for these instances were randomly generated within the interval [1,10][1,10]. The coordinate ranges for the points in these instances are provided in the third column of Table \ref{inst-Ruspini-Bongartz}.

The results for all test problems are reported using the squared Euclidean norm with a stopping criterion of $||\bar{\mathbf{{x}}}-\mathbf{x}^{(k)}||\leq0.0001$, as well as the Euclidean and rectilinear norms with $\dfrac{|\bar{F}-F_k|}{\bar{F}}\leq0.01$. The final iteration is denoted by $k$.

\begin{table}\caption{ The coordinate ranges of existing points of instances from \cite{B90}.}\label{inst-Ruspini-Bongartz}
	\centering
	\small 
	\begin{tabular}{ccc}
		\hline
		Instance&$n$&Coordinate ranges\\
		\hline
		Ruspini&75&$(4,4)\ to\ (117,156)$\\
		\hline
		Bongartz&287&$(5,5)$ to $(48,48)$\\
		\hline
	\end{tabular}
\end{table}


The results of testing our algorithms on the instances using the squared Euclidean norm are presented in Tables \ref{result-SE-18data-(2,2)} to \ref{result-SE-Ruspini-Bongartz-CPU}, where they are compared with those obtained by the linear programming method of Baroughi-Bonab et al. \cite{BBA10}. In these tables, $\mathbf{x}^{(k)}$ and $C^*$ represent the final point and optimal cost achieved by our algorithms, respectively. Additionally, $C_B$ denotes the optimal cost obtained by the method of Baroughi-Bonab et al. \cite{BBA10}, and $\mathbf{{x}}_B$ is the optimal solution derived by applying the center of gravity method to the coordinates generated by their approach.

A comparison of the results reveals that the total cost produced by our algorithms is approximately equal to that of the Baroughi-Bonab et al. method \cite{BBA10}. However, their method is faster than our algorithms, and the solution $\mathbf{{x}}_{B}$ obtained by their approach is exact,  i.e., $||\bar{\mathbf{{x}}}-\mathbf{{x}}_{B}||=0.0000$. It is important to note that their method is specifically designed for the squared Euclidean norm and cannot be applied to other norms.

 
\begin{table}\caption{ The obtained coordinates for the instance with 18 points with squared Euclidean norm for the case $\bar{\mathbf{x}}=(2,2)$.}		\label{result-SE-18data-(2,2)}
	\centering
	\small 
	\begin{tabular}{|c||c|c|}
		\hline
		Algorithm&ISFLP1 and ISFLP2&Baroughi et al. \cite{BBA10}\\
		\hline
		$\mathbf{x}^{(K)}$&$(2.0001 , 2.0000)$&-\\
		\hline
		$C^*$&$78.3324$&-\\
		\hline
		${\mathbf{x}_B}$&-&$(2,2)$\\
		\hline
		$C_B$&-&$78.3333$\\
		\hline
		k&$21$&-\\
		\hline
		$\mathbf{A}_{1}^{(k)}=(a^{(k)}_1,b^{(k)}_1)$&
		(-42.6660 , 2.0000)&(-42.6667 , 2.0000)\\
		$\mathbf{A}_{2}^{(k)}=(a^{(k)}_2,b^{(k)}_2)$&
		(1.0000 , 3.0000)&(1.0000 , 3.0000)\\
		$\mathbf{A}_{3}^{(k)}=(a^{(k)}_3,b^{(k)}_3)$&
		(2.0000 , 5.0000)&(2.0000 , 5.0000)\\
		$\mathbf{A}_{4}^{(k)}=(a^{(k)}_4,b^{(k)}_4)$&
		(3.0000 , 6.0000)&(3.0000 , 6.0000)\\
		$\mathbf{A}_{5}^{(k)}=(a^{(k)}_5,b^{(k)}_5)$&
		(4.0000 , 8.0000)&(4.0000 , 8.0000)\\
		$\mathbf{A}_{6}^{(k)}=(a^{(k)}_6,b^{(k)}_6)$&
		(4.0000 , 1.0000)&(4.0000 , 1.0000)\\
		$\mathbf{A}_{7}^{(k)}=(a^{(k)}_7,b^{(k)}_7)$&
		(4.0000 , 4.0000)&(4.0000 , 4.0000)\\
		$\mathbf{A}_{8}^{(k)}=(a^{(k)}_8,b^{(k)}_8)$&
		(4.0000 , 9.0000)&(4.0000 , 9.0000)\\
		$\mathbf{A}_{9}^{(k)}=(a^{(k)}_9,b^{(k)}_9)$&
		(5.0000 , 3.0000)&(5.0000 , 3.0000)\\
		$\mathbf{A}_{10}^{(k)}=(a^{(k)}_{10},b^{(k)}_{10})$&(5.0000 , 5.0000)&(5.0000 , 5.0000)\\
		$\mathbf{A}_{11}^{(k)}=(a^{(k)}_{11},b^{(k)}_{11})$&(6.0000 , 6.0000)&(6.0000 , 6.0000)\\
		$\mathbf{A}_{12}^{(k)}=(a^{(k)}_{12},b^{(k)}_{12})$&(6.0000 , 3.0000)&(6.0000 , 3.0000)\\
		$\mathbf{A}_{13}^{(k)}=(a^{(k)}_{13},b^{(k)}_{13})$&(7.0000 , 1.0000)&(7.0000 , 1.0000)\\
		$\mathbf{A}_{14}^{(k)}=(a^{(k)}_{14},b^{(k)}_{14})$&(7.0000 , 2.0000)&(7.0000 , 2.0000)\\
		$\mathbf{A}_{15}^{(k)}=(a^{(k)}_{15},b^{(k)}_{15})$&(8.0000 , 5.0000)&(8.0000 , 5.0000)\\
		$\mathbf{A}_{16}^{(k)}=(a^{(k)}_{16},b^{(k)}_{16})$&(8.0000 , 8.0000)&(8.0000 , 8.0000)\\
		$\mathbf{A}_{17}^{(k)}=(a^{(k)}_{17},b^{(k)}_{17})$&(9.0000 , -27.6664)&(9.0000 , -26.2500)\\
		$\mathbf{A}_{18}^{(k)}=(a^{(k)}_{18},b^{(k)}_{18})$&(9.0000 , 6.0000)&(9.0000 , 6.0000)\\
		\hline
		$||\bar{\mathbf{{x}}}-\mathbf{x}^{(k)}||$&$0.00005$&-\\
		\hline
		$||\bar{\mathbf{{x}}}-\mathbf{x}_B||$&-&$0.00000$\\
		\hline
		$CPU( in \ \ sec)$&$0.6667$ and $0.4800$&$0.0821$\\
		\hline
		\end{tabular}
\end{table}

\begin{table}\caption{ The results 
		for the instances with squared Euclidean norm.}\label{result-SE-18data}
	\centering
	\tiny
	\begin{tabular}{ccc|ccc|ccc|cc}
		\hline
 		&&&\multicolumn{3}{c|}{ISFLP1}&\multicolumn{3}{c|}{ISFLP2}&\multicolumn{2}{c}{Method of Baroughi et al. \cite{BBA10}}\\ 
		\hline
		Instance&$n$&$\mathbf{\bar{x}}=(\bar{x},\bar{y})$&$C^*$&$iter$&$||\bar{\mathbf{{x}}}-\mathbf{x}^{(k)}||$&$C^*$&$iter$&$||\bar{\mathbf{{x}}}-\mathbf{x}^{(k)}||$&$C_B$&$||\bar{\mathbf{{x}}}-\mathbf{x}_B||$\\
		\hline
		data &18&$(7,7)$&$54.9989$&$20$
		&$0.00005$&$54.9989$&$20$&$0.00005$&$55$&$0.00000$\\
		& &$(-3,-5)$&$238.3323$&$21$&$0.00007$&
		$238.3323$&$21$&$0.00007$&$238.3333$&$0.00000$\\
	\hline
	Ruspini&75&$(50,50)$&$1.35e+03$&$23$&$0.00007$&$5.93e+04$&$46$&$0.00006$&$1.35e+03$&$0.00000$\\
	& &$(-80,-20)$&$6.13e+03$&$26$&$0.00008$&$6.13e+03$&$26$&$0.00008$&$6.13e+03$&$0.0000$\\
	& &$(-20,80)$&$2.07e+03$&$24$&$0.00006$&$2.07e+03$&$24$&$0.00006$&$2.06e+03$&$0.00000$\\
	\hline
	Bongartz&287&$(15,35)$&$8.75e+02$&$21$&$0.00006$&$8.75e+02$&$21$&$0.00008$&$8.73e+02$&$0.00000$\\
	& &$(20,29)$&$4.17e+02$&$20$&$0.00009$&$4.17e+02$&$20$&$0.00005$&$4.16e+02$&$0.00000$\\
	& &$(50,45)$&$3.90e+03$&$24$&$0.00009$&$3.90e+03$&$24$&$0.00008$&$3.88e+03$&$0.00000$\\
		\hline
	\end{tabular}
\end{table}


\begin{table}\caption{ Average CPU(in sec) times of results in Table \ref{result-SE-18data}.}\label{result-SE-Ruspini-Bongartz-CPU}
	\centering
	\tiny
	\begin{tabular}{cc|c|c|c}
		\hline
		Instance&$n$&\multicolumn{1}{c|}{ISFLP1}&\multicolumn{1}{c|}{ISFLP2}&\multicolumn{1}{c}{Method of Baroughi et al. \cite{BGG10}}\\ 
		\hline
		Ruspini&75&1.95&$1.36$&$0.05$\\
		\hline
		Bongartz&287&2.45&$1.64$&$0.05$\\
		\hline
	\end{tabular}
\end{table}


Tables \ref{result-E-18data-(2,2)} to \ref{result-L1-18data} contain the results obtained by the Algorithms ISFLP1 and ISFLP2 for Euclidean and rectilinear norms and given points $\mathbf{\bar{x}}$.

\begin{table}\caption{ The coordinates obtained by ISFLP1 and ISFLP2 for the instance with 18 points with Euclidean norm in the case $\bar{\mathbf{x}}=(2,2)$.}\label{result-E-18data-(2,2)}
	\centering
	\small 
	\begin{tabular}{|c||c|c|}
		\hline
		$Algorithm$&$ISFLP1$&$ISFLP2$\\
		\hline
		$C^*$&80.8838&80.7814\\
		\hline
		k&3&3\\
		\hline
		$\mathbf{A}_{1}^{(k)}=(a^{(k)}_1,b^{(k)}_1)$
		&(1.0647 , 1.7363)&(1.0645 , 1.7372)\\
		$\mathbf{A}_{2}^{(k)}=(a^{(k)}_2,b^{(k)}_2)$
		&(0.9775 , 2.4035)&(0.9773 , 2.4060)\\
		$\mathbf{A}_{3}^{(k)}=(a^{(k)}_3,b^{(k)}_3)$
		&(1.9631 , 5.0009)&(1.9632 , 5.0009)\\
		$\mathbf{A}_{4}^{(k)}=(a^{(k)}_4,b^{(k)}_4)$
		&(2.0000 , 2.0000)&(2.0000 , 2.0000)\\
		$\mathbf{A}_{5}^{(k)}=(a^{(k)}_5,b^{(k)}_5)$
		&(3.8981 , 7.9291)&(3.8986 , 7.9291)\\
		$\mathbf{A}_{6}^{(k)}=(a^{(k)}_6,b^{(k)}_6)$
		&(1.7250 , 1.1160)&(1.7262 , 1.1147)\\
		$\mathbf{A}_{7}^{(k)}=(a^{(k)}_7,b^{(k)}_7)$
		&(3.9533 , 4.1297)&(3.9534 , 4.1294)\\
		$\mathbf{A}_{8}^{(k)}=(a^{(k)}_8,b^{(k)}_8)$
		&(1.9218 , 1.9568)&(1.9218 , 1.9565)\\
		$\mathbf{A}_{9}^{(k)}=(a^{(k)}_9,b^{(k)}_9)$
		&(1.9986 , 1.9991)&(1.9978 , 1.9985)\\
		$\mathbf{A}_{10}^{(k)}=(a^{(k)}_{10},b^{(k)}_{10})$
		&(4.9612 , 5.2475)&(4.9612 , 5.2469)\\
		$\mathbf{A}_{11}^{(k)}=(a^{(k)}_{11},b^{(k)}_{11})$
		&(5.9934 , 6.1629)&(5.9934 , 6.1626)\\
		$\mathbf{A}_{12}^{(k)}=(a^{(k)}_{12},b^{(k)}_{12})$
		&(5.8310 , -2.2198)&(5.8317 , -2.2108)\\
		$\mathbf{A}_{13}^{(k)}=(a^{(k)}_{13},b^{(k)}_{13})$
		&(1.8947 , 1.8861)&(1.8941 , 1.8861)\\
		$\mathbf{A}_{14}^{(k)}=(a^{(k)}_{14},b^{(k)}_{14})$
		&(6.9350 , -0.3228)&(6.9354 , -0.3034)\\
		$\mathbf{A}_{15}^{(k)}=(a^{(k)}_{15},b^{(k)}_{15})$
		&(8.0635 , 5.0219)&(8.0635 , 5.0219)\\
		$\mathbf{A}_{16}^{(k)}=(a^{(k)}_{16},b^{(k)}_{16})$
		&(7.9739 , 8.0306)&(7.9739 , 8.0306)\\
		$\mathbf{A}_{17}^{(k)}=(a^{(k)}_{17},b^{(k)}_{17})$
		&(9.1629 , 6.8984)&(9.1630 , 6.8983)\\
		$\mathbf{A}_{18}^{(k)}=(a^{(k)}_{18},b^{(k)}_{18})$
		&(9.0404 , 5.9794)&(9.0404 , 5.9794)\\
		\hline
		$\frac{|\bar{F}-F_k|}{|\bar{F}|}$
		&0.0070&0.0071\\
		\hline
		$CPU( in \ \ sec)$&35.4190&31.4453\\
		\hline
	\end{tabular}
\end{table}

\begin{table}\caption{ The results of 
		for the instances with Euclidean norm.}\label{result-E-18data}
	\centering
	\tiny
	\begin{tabular}{ccc|cccc|cccc}
		\hline
		&&&\multicolumn{4}{c|}{ISFLP1}&\multicolumn{4}{c}{ISFLP2}\\ 
		\hline
		
		Instance&$n$&$\mathbf{\bar{x}}=(\bar{x},\bar{y})$&$C^*$&$iter$&$\dfrac{|\bar{F}-F_k|}{|\bar{F}|}$&$CPU$&$C^*$&$iter$&$\dfrac{|\bar{F}-F_k|}{|\bar{F}|}$&$CPU$\\
		&&&&&&(in sec)&&&&(in sec)\\
		\hline
		data &18&$(7,7)$&$22.8798$&$3$
		&$0.0080$&$34.9575$&$22.6452$&3&0.0088&$32.0751$\\
		& &$(-3,-5)$&$390.2082$&$4$&$0.0034$&$55.3343$
		&$360.2234$&$4$&$0.0037$&$46.2091$\\
		\hline
		Ruspini&75&$(50,50)$&2.77e+03&4&$0.0048$&$161.72$&2.65e+03&4&0.0074&121.22\\
		& &$(-80,-20)$&3.25e+04&5&0.0032&$220.70$&3.30e+04&5&0.0032&151.66\\
		& &$(-20,80)$&7.10e+03&5&0.0050&$230.35$&7.28e+03&5&0.0040&169.99\\
		\hline
		Bongartz&287&$(15,35)$&4.04e+03&12&0.0084&$3033.47$&4.08e+03&13&0.0071&1278.28\\
		& &$(20,29)$&1.17e+03&6&0.0089&748.78&1.27e+03&7&0.0077&623.36\\
		& &$(50,45)$&2.17e+04&36&0.0082&79334.95&2.18e+04&43&0.0079&32920.78\\
		\hline
	\end{tabular}
\end{table}


\begin{table}\caption{ The obtained coordinates 
		for the instance of 18 points with rectilinear norm in the case $\bar{\mathbf{x}}=(2,2)$. }\label{result-L1-18data-(2,2)}
	\centering
	\small 
	\begin{tabular}{|c||c|c|}
		\hline
		$Algorithm$&$ISFLP1$&$ISFLP2$\\
		\hline
		\hline
		$C^*$&$85.7698$&$85.5271$\\
		\hline
		k&$5$&$5$\\
		\hline
		$\mathbf{A}_{1}^{(k)}=(a^{(k)}_1,b^{(k)}_1)$&(0.7963 , 1.9825)&(0.8595 , 1.9825)\\
		$\mathbf{A}_{2}^{(k)}=(a^{(k)}_2,b^{(k)}_2)$&(1.0000 , 1.9889)&(1.0000 ,1.9975)\\
		$\mathbf{A}_{3}^{(k)}=(a^{(k)}_3,b^{(k)}_3)$&(1.9819 , 4.9462)&(1.9858 , 4.9376)\\
		$\mathbf{A}_{4}^{(k)}=(a^{(k)}_4,b^{(k)}_4)$&(1.9988 , 6.2974)&(1.9864 , 6.2060)\\
		$\mathbf{A}_{5}^{(k)}=(a^{(k)}_5,b^{(k)}_5)$&(2.2953 , 7.7026)&(1.9784 , 7.7940)\\
		$\mathbf{A}_{6}^{(k)}=(a^{(k)}_6,b^{(k)}_6)$&(1.9887 , 1.0415)&(1.9887 , 1.0330)\\
		$\mathbf{A}_{7}^{(k)}=(a^{(k)}_7,b^{(k)}_7)$&(1.9990 , 4.1863)&(2.0000 , 4.1279)\\
		$\mathbf{A}_{8}^{(k)}=(a^{(k)}_8,b^{(k)}_8)$&(1.7001 , 8.6238)&(1.8043 , 8.7344)\\
		$\mathbf{A}_{9}^{(k)}=(a^{(k)}_9,b^{(k)}_9)$&(1.9974 , 1.9441)&(1.9854 , 1.9556 )\\
		$\mathbf{A}_{10}^{(k)}=(a^{(k)}_{10},b^{(k)}_{10})$&(1.9896 , 1.9573)&(1.9898 , 1.9749)\\
		$\mathbf{A}_{11}^{(k)}=(a^{(k)}_{11},b^{(k)}_{11})$&(6.0592 , 6.1561)&(6.0477 , 6.1215)\\
		$\mathbf{A}_{12}^{(k)}=(a^{(k)}_{12},b^{(k)}_{12})$&(6.0094 , 1.9468)&(5.9983 , 1.9546)\\
		$\mathbf{A}_{13}^{(k)}=(a^{(k)}_{13},b^{(k)}_{13})$
		&(6.8854 , 1.1146)&(6.9115 , 1.0885)\\
		$\mathbf{A}_{14}^{(k)}=(a^{(k)}_{14},b^{(k)}_{14})$
		&(7.0358 , 1.9804)&(7.0289 , 1.9629)\\
		$\mathbf{A}_{15}^{(k)}=(a^{(k)}_{15},b^{(k)}_{15})$
		&(8.1381 , 2.0000)&(8.1073 , 2.0000)\\
		$\mathbf{A}_{16}^{(k)}=(a^{(k)}_{16},b^{(k)}_{16})$
		&(8.0000 , 8.0565)&(8.0000 , 8.0445)\\
		$\mathbf{A}_{17}^{(k)}=(a^{(k)}_{17},b^{(k)}_{17})$
		&(9.4815 , 6.6534)&(9.3488 , 6.7575)\\
		$\mathbf{A}_{18}^{(k)}=(a^{(k)}_{18},b^{(k)}_{18})$
		&(9.0937 , 5.9642)&(9.0655 , 5.9711)\\
		\hline
		$\dfrac{|\bar{F}-F_k|}{|\bar{F}|}$
		&0.0054&0.0000\\
		\hline
		$CPU( in \ \ sec)$&$22.1654$&$55.3114$\\
		\hline
		\end{tabular}
\end{table}

\begin{table}\caption{ The results 
		for the instances 
		with rectilinear norm.}\label{result-L1-18data}
	\centering
	\tiny
	\begin{tabular}{ccc|cccc|cccc}
		\hline
		&&&\multicolumn{4}{c|}{ISFLP1}&\multicolumn{4}{c}{ISFLP2}\\ 
		\hline
		Instance&$n$&$\mathbf{\bar{x}}=(\bar{x},\bar{y})$&$C^*$&$iter$&$\dfrac{|\bar{F}-F_k|}{|\bar{F}|}$&CPU&$C^*$&$iter$&$\dfrac{|\bar{F}-F_k|}{|\bar{F}|}$&CPU\\
		&&&&&&(in sec)&&&&(in sec)\\
		\hline
		data &18&$(7,7)$&$24.4908$&$3$
		&$0.0002$&13.6394&$24.4702$&$3$&$0.0003$&$35.8552$\\
		& &$(-3,-5)$&$484.1516$&4&0.0000
		&$14.4387$&483.7831&4&0.0000&32.8848\\
		\hline
		Ruspini&75&$(50,50)$&2.27e+03&4&0.0034&35.57&2.03e+03&3&0.0088&58.98\\
		& &$(-80,-20)$&4.67e+04&4&0.0012&35.57&4.56e+04&4&0.0026&85.16\\
		& &$(-20,80)$&1.04e+04&3&0.0076&25.91&1.06e+04&3&0.0070&58.60\\
		\hline
		Bongartz&287&$(15,35)$&3.21e+03&13&0.0082&1585.68&3.16e+03&13&0.0085&1973.37\\
		& &$(20,29)$&1.12e+03&8&0.0098&579.71&1.11e+03&8&0.0004&830.93\\
		& &$(50,45)$&2.32e+04&37&0.0093&57600.85&2.31e+04&41&0.0084&42655.04\\
		\hline
		\end{tabular}
\end{table}


\subsection{Analysis  on the effect of given point}
In this section, we examine the effect of given point $\mathbf{\bar{x}}$ on the number of iterations of Algorithms ISFLP1 and ISFLP2. Obviously, if  $\mathbf{\bar{x}}$ is the optimal solution of the considered location problem with initial given coordinates, then algorithms stop in the first iteration. In inverse minisum location problem with variable weights, if the given point $\mathbf{\bar{x}}$ does not lie in the convex hull of  $\mathbf{P}_1 , ... , \mathbf{P}_n$, then the inverse model is infeasible \cite{F22}. However, this does not apply to inverse minisum location problem with variable coordinates. As the reported result in tables \ref{result-SE-4points} to \ref{result-L1-18data} show, there is not any meaningful relation between choosing given point $\mathbf{\bar{x}}$ and the number of iterations of the algorithms. For example, consider instance with 18 points in table \ref{inst-data-18}. The optimal solution of this problem with the given coordinates for squared Euclidean and rectilinear norms are points $\mathbf{x}^{(0)}=(5.2750 , 4.6000)$ and $\mathbf{x}^{(0)}=(5.0000 , 5.0000)$, respectively. The results indicate that if $\mathbf{\bar{x}}=(2 , 2)$ and $\mathbf{\bar{x}}=(-3 ,-5)$ then Algorithms ISFLP1 and ISFLP2 terminate after 21 iterations with squared Euclidean norm (see Tables \ref{result-SE-18data-(2,2)} and \ref{result-SE-18data}). Also, algorithms terminate after 5 and 4 iterations with $L_1$ norm, respectively (see Tables \ref{result-L1-18data-(2,2)} and \ref{result-L1-18data}). Whereas obviously, the first point is near than the second one to the optimal solution. Similarly, consider Ruspini instance with $L_2$ norm and optimal solution of $\mathbf{x}^{(0)}=(50.2408 , 103.4575)$. The results indicate that if $\mathbf{\bar{x}}=(-80 , -20)$ and $\mathbf{\bar{x}}=(-20 , 80)$ then Algorithms ISFLP1 and ISFLP2 terminate after 5 iterations.
    
The numerical results indicate that the number of iterations in Algorithm ISFLP1 is less than or equal to that of Algorithm ISFLP2. Algorithm ISFLP2 has a lower time complexity than Algorithm ISFLP1. However, we believe that Algorithm ISFLP1 performs better than Algorithm ISFLP2. Consider Example \ref{ex1} with $L_2$ norm and $\mathbf{\bar{x}}=(-2 , 5)$. The numerical results in Table \ref{result-E-4points-(-2,5)} illustrate this point.

    \begin{table}\caption{ The results of Algorithms ISFLP1 and ISFLP2
    		for the instances from Example \ref{ex1} with Euclidean norm.}\label{result-E-4points-(-2,5)}
    	\centering
    	\tiny
    	\begin{tabular}{ccc|ccc|ccc}
    		\hline
    		&&&\multicolumn{3}{c|}{ISFLP1}&\multicolumn{3}{c}{ISFLP2}\\ 
    		\hline
    		
    		Instance&$n$&$\mathbf{\bar{x}}=(\bar{x},\bar{y})$&$C^*$&$iter$&CPU&$C^*$&$iter$&CPU\\
    		&&&&&(in sec)&&&(in sec)\\
    		\hline
    		Example \ref{ex1} &4&$(-2,5)$&$39.3674$&$39$
    	&273.8595&$39.1801$&$2111$&$6181.4753$\\
    	 		\hline
    		\end{tabular}
    \end{table}

\subsection{Experimental notes }
In summary, considering the numerical experiments and discussions in the previous sections, the following notes can be stated:
\begin{enumerate}
\item  Any of three stopping conditions  $|f(\mathbf{\bar{x}},\mathbf{w},\mathbf{P}^{(t)})-f(\mathbf{{x}}^{(t)},\mathbf{w},\mathbf{P}^{(t)})|\leq \epsilon$, $||\mathbf{\bar{x}}-\mathbf{{x}}^{(t)}||\leq \epsilon$, and  $||\mathbf{P}^{(t)}-\mathbf{P}^{(t-1)}||\leq \epsilon$ can be used for the minisum case.
\item If $\mathbf{\bar{x}}$ is not the optimal, then its location does not affect on the number of iterations of the algorithm.
\item  Figures \ref{Initial-Points-Bongartz} to \ref{Final-Bongartz-L1} show the initial coordinates and results of Algorithm ISFLP1 with Euclidean and rectilinear norms for Bongartz 287 instance and $\bar{\mathbf{x}}=(50,45)$ , respectively. Figures show that Algorithm ISFLP1 results in optimal coordinates around the given point $\bar{\mathbf{x}}$.
\begin{figure}
	\centering
	\includegraphics[scale=0.6]{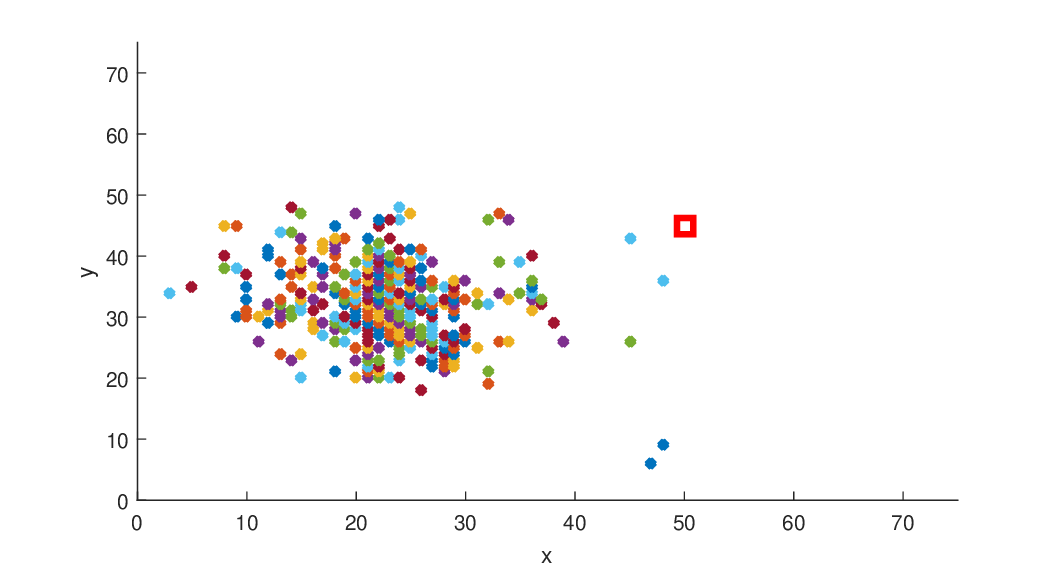}
	\caption{Initial coordinates of Bongartz and $\bar{\mathbf{x}}=(50,45)$}\label{Initial-Points-Bongartz} 
\end{figure}

\begin{figure}
	\centering
	\includegraphics[scale=0.6]{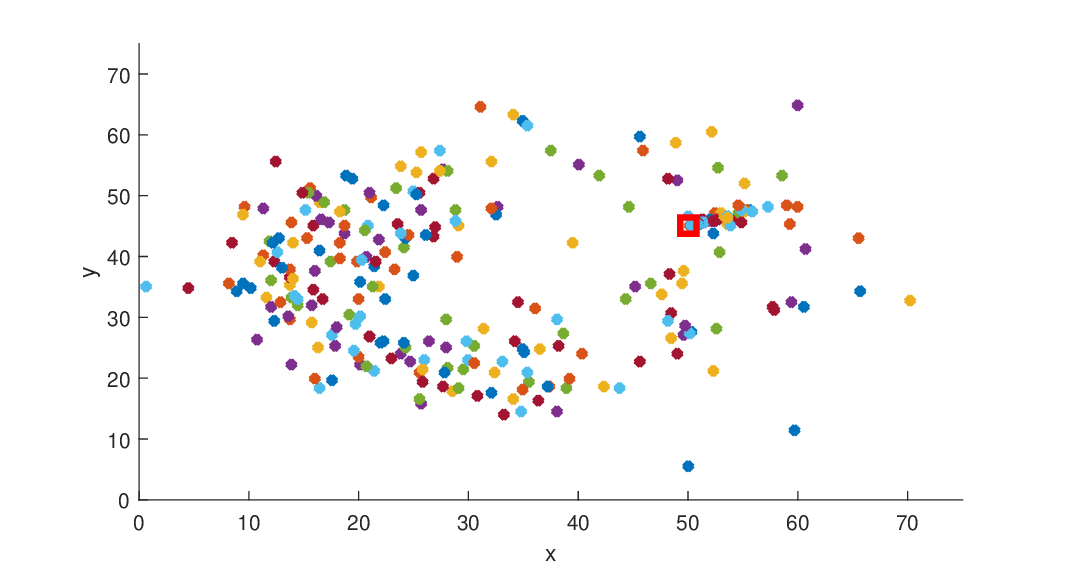}
	\caption{The results of Algorithms ISFLP1 for Bongartz 287 instance with Euclidean norm and $\bar{\mathbf{x}}=(50,45)$} \label{Final-Bongartz-l2}
\end{figure}

\begin{figure}
	\centering
	\includegraphics[scale=0.6]{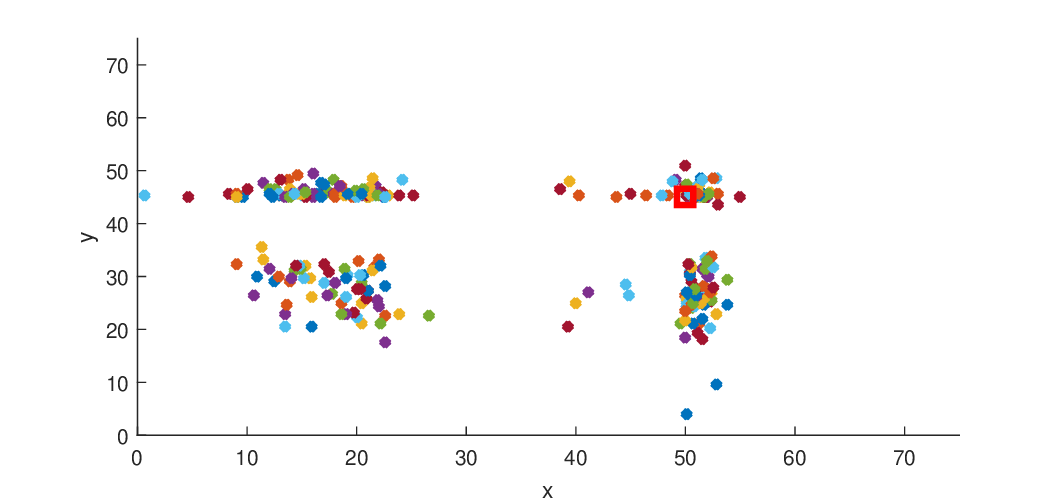}
	\caption{The results of Algorithms ISFLP1 for Bongartz 287 instance with rectilinear norm and $\bar{\mathbf{x}}=(50,45)$} \label{Final-Bongartz-L1}
\end{figure}
\end{enumerate}  	

\section{Summary and conclusion}\label{sec5}

In this paper we presented two algorithms for general case of inverse continuous location problems. The optimality condition have been proved for the given algorithm. In special cases we considered the inverse minisum
 location models and some test problems have been solved by the presented algorithms. The results show the algorithms could efficiently find the optimal solutions.
 
Note that the presented algorithms can be applied for solving inverse of other cases of location problems such as the inverse covering problem, inverse line location problem and inverse circle location problem. 



\end{document}